\documentclass[12pt]{article}

\usepackage{amssymb}
\usepackage[dvips]{graphics}

\newcommand{\vv}{\vspace{2mm}}

\def\Z{{\bf Z}}

\def\Aut{{\mathrm{Aut}}}

\def\pf{{\indent\textit{Proof.}\ }}
\def\qed{\hfill$\square$}
\newcounter{para}[section]
\renewcommand{\thepara}{\thesection.\arabic{para}}
\renewcommand{\thesection}{\arabic{section}}
\renewcommand{\paragraph}{\refstepcounter{para}
\indent{\bf{\thepara}}}
\newcommand{\sectioning}{\refstepcounter{section}
\indent{\bf \thesection.}}

\newcommand{\an}{{\mbox{{\rm {\tiny an}}}}}

\begin{document}

{\footnotesize \noindent Running head: Groups of Lam\'e type \hfill July 26, 2002

\noindent Math.\ Subj.\ Class.\ (2000): 14G22, 14H37, 20E08}

\vspace{40mm}

\begin{center} 

{\Large Mumford curves with maximal automorphism group II:

\vv

Lam\'e type groups in genus 5-8}

\vv

{\sl by} Gunther Cornelissen {\sl and} Fumiharu Kato

\vv

\end{center}

\vv

\vv

{\bf Introduction} 

\vv

\vv

\noindent  Let $(k,|\cdot|)$ be a non-archimedean valued field of characteristic $p \geq 0$, and $X$ a Mumford curve of genus $g$ over $k$ (i.e., such that its stable reduction is a union of rational curves intersecting in rational points over the residue class field $\bar{k}$ of $k$). It is known (cf.\ \cite{Herrlich:80}, \cite{CKK}) that for $g \in \{ 5,6,7,8 \}$, $X$ has at most $12(g-1)$ automorphisms (where we assume for $p=0$ that char$(\bar{k})>5$). The aim of this work is to determine 
explicitly all (families of) Mumford curves that attain this bound. The main 
result is the following:

\vv

{\bf Theorem.} \  {\sl Let $X$ be a Mumford curve of genus $g \in \{5,6,7,8\}$ with $12(g-1)$ automorphisms over a non-archimedean field $(k,|\cdot |)$, such that char$(\bar{k})>5$ if $p=0$ Then actually $p>3$ and either 

}(a) {\sl $g=5$ and $\Aut(X) \cong S_4 \times {\bf Z}_2$; or

}(b) {\sl $g=6$, $\Aut(X) \cong A_5$ and $p \neq 5$. 

\noindent Furthermore, the normalizers of the corresponding Schottky groups are isomorphic to the following amalgams of groups: 

in case }(a) {\sl to either }(a1) {\sl $S_4 *_{{\bf Z}_4} D_4$ or }(a2) {\sl $D_3 *_{{\bf Z}_2} D_2$;

in case }(b) {\sl to either
}(b1) {\sl $A_5 *_{{\bf Z}_5} D_5$, }(b2) {\sl  $A_4 *_{{\bf Z}_3} D_3$ or 
}(b3) {\sl $D_3 *_{{\bf Z}_2} D_2$. }

\vv

We then turn to determining the strata in the moduli space of curves of genus $g$ with automorphism group as above, and then determine the loci of Mumford
curves in these strata.

\vv

For case (a), we have the following result:  

\vv 

{\bf Proposition A.} \ {\sl Let $X$ be a Mumford curve of genus 5 with automorphism group $S_4 \times \Z_2$. Then $X$ occurs in the family of curves $C_\alpha$ constructed as follows. Let $C^\sim_\alpha \rightarrow {\bf P}^1_\alpha$ be the following one-parameter family of genus 3 curves $$ C^\sim_\alpha \ : \ x^4+y^4+z^4+\alpha(x^2y^2+y^2z^2+x^2z^2)=0; $$
then $\Aut(C^\sim_\alpha)$ contains $S_4$, and 
the 4 lines $ \ x \pm y \pm z = 0$ are bitangent to $C^\sim_\alpha$ at 8 points, which form two disjoint $S_4$-orbits, say, $\{p_i\}_{i=1}^4$ and $\{p'_i\}_{i=1}^4$. The divisor $\sum_{i=1}^4 (p_i - p'_i)$ corresponds to an $S_4$-invariant two-torsion point of $\mathrm{Jac}(C^\sim_\alpha)$, so defines an \'etale ${\bf Z}_2$-cover $C_\alpha$ of $C^\sim_\alpha$. 

More precisely, $C_\alpha$ is
a Mumford curve exactly for 

\textup{(a1)} $|\alpha|>1$ or \ \textup{(a2)} $|\alpha+2|<1$, 

\noindent where the numeration is compatible with that of the normalizers in the theorem.}

\vv

The curves in (b) belong to the one-parameter family of genus six curves on the Del Pezzo quintic that was studied by Edge in \cite{Edge1} \cite{Edge2}, and its Mumford loci were considered in  \cite{Kato:01}. For the sake of completeness, 
we state the result.

\vv

{\bf Proposition B.} \ {\sl Let $X$ be a Mumford curve of genus $g=6$ with automorphism group $A_5$. Then $X$ occurs in the family $E_\alpha$ of genus 6 curves in the Del Pezzo quintic surface which is the strict transform of the family of sextics $E^\sim_\alpha$ in ${\bf P}^2$ given 
by $E^\sim_\alpha  :  T + \alpha S = 0$ where
\begin{eqnarray*} & &  T= x^6+y^6+z^6+(x^2+y^2+z^2)(x^4+y^4+z^4)-12x^2y^2z^2 \\
& &  S= (y^2-z^2)(z^2-x^2)(x^2-y^2)
\end{eqnarray*}
More precisely, $E_\alpha$ is a Mumford curve precisely when $\alpha$ is in one of the loci 

\textup{(b1)} $0<|\alpha \pm 5 \sqrt{5}| < 1$; \ \textup{(b2)} $0<|\alpha \pm \sqrt{-3}|<1$ or \ \textup{(b3)} $|\alpha|>1$

\noindent where the numeration is compatible with that of the normalizers in the theorem (and the two different signs in \textup{(b1)} and \textup{(b2)} actually correspond to the same locus in moduli space).}

\vv

\vv

Notice that for $g \notin \{5,6,7,8\}$ and $p>0$, the maximal number of automorphisms of a Mumford curve of genus $g \geq 2$ is $2 \sqrt{g} (\sqrt{g}+1)^2$ (cf.\ \cite{CKK}), and in \cite{CK}, the corresponding families of curves that attain this bound were explicitly found. Thus, the above theorem and proposition complete the determination of all positive characteristic Mumford curves with maximal automorphism group begun in \cite{CK}, i.e., the non-archimedean analogue of the determination of all Hurwitz groups (cf.\ Conder, \cite{Conder}). 

Here is a short outline of the proofs. We first show that $X \mapsto \Aut(X) \backslash X$ is a cover of ${\bf P}^1$ ramified tamely above four points with indices $(2,2,2,3)$ (this follows from manipulating the Hurwitz formula for ordinary curves and uses the main result from \cite{CKK}).  Once the ramification type has been fixed, we use the techniques of \cite{CKK} to make a finite
list of possible (abstract types of) normalizers of the corresponding Schottky groups (there turn out to be only four). Namely, such $N$ form a tree product,
which we can combinatorially rewrite as a product over an easier tree with only four ends (corresponding to the four ramification points).
We then look at what finite groups of order $12(g-1)$ can be quotients of these $N$. Actually, most can be excluded using elementary, but long-winded group theory (of which we defer some to an appendix). To explicitly determine the corresponding $S_4 \times \Z_2$-family, we decompose it into an \'etale $\Z_2$-part and a genus three $S_4$-cover of ${\bf P}^1$. We then use classical algebraic geometry to make the cover explicit as in the proposition, and then determine the Mumford loci in the moduli by computing a semistable model for the degenerate fibers by the methods from \cite{Kato:01}.

Notice that a Riemann surface has at most $84(g-1)$ automorphisms, and if
this occurs, then it corresponds to a $(2,3,7)$-cover of ${\bf P}^1$ --- possible automorphism groups are then called {\sl Hurwitz groups}. Here, we
are dealing with $(2,2,2,3)$-covers of ${\bf P}^1$ by ``non-archimedean Riemann surfaces'', and as these ramification indices remind us of the Lam\'e equation, we would like to baptize the corresponding possible automorphism groups {\sl of Lam\'e type}.

\vv

\vv

\sectioning\label{prel}\ 
{\bf Classification of normalizers of Schottky groups}

\vv

\vspace{1ex}
\paragraph\label{Mum} {\bf Mumford curves} (\cite{Mumford:72}). \ Let $(k,|\cdot|)$ be a non-archimedean valued field of characteristic $p \geq 0$, and $X$ a Mumford curve
of genus $g$ over $k$. This means that there exists a semi-stable model over 
the valuation ring of $k$ such that
the stable reduction of $X$ over the residue field $\bar{k}$ of $k$ is a union of rational curves intersecting 
in $\bar{k}$-rational points. Equivalently, as a rigid analytic space over $k$, the analytification
$X^{\an}$ of $X$ is isomorphic to an analytic space of the form $\Gamma \backslash ({\bf P}^{1,\an}_k-\mathcal{L})$,
where $\Gamma$ (the so-called Schottky group of $X$) is a discrete free subgroup of $PGL(2,k)$ of rank 
$g$ (acting in the obvious way on ${\bf P}^{1,\an}_k$) with $\mathcal{L}$ as its set of limit points.

\vspace{1ex}
\paragraph\label{RH} {\bf Proposition.} \  {\sl Let $X$ be a Mumford curve of genus $g \in \{5,6,7,8\}$ over a non-archimedean field $(k,|\cdot |)$ (with $\mbox{{\rm char}}(\bar{k})>5$ if $p=0$), such that the order of its automorphism group $\Aut(X)$ is $12(g-1)$.  Then $p=0$ or $p>3$,  the quotient curve $\Aut(X) \backslash X$ is isomorphic to ${\bf P}^1$ and the ramification in the corresponding quotient map is tame of type $(2,2,2,3)$.}

\vv

\pf It follows from the proof of Satz 3 in \cite{Stichtenoth} that if $|\Aut(X)|=12(g-1)$, then $Y:=\Aut(X) \backslash X$ if of genus zero and 
$X \rightarrow Y$ ramifies above at most four points. However, if it branches above strictly less than four points, then it is known for $p>0$ that $|\Aut(X)| \leq 2 \sqrt{g} (\sqrt{g}+1)^2$ (cf.\ \cite{CKK}, \S 6), which is stricly less than $12(g-1)$ for $g$ in the prescribed range. On the other hand, if $p=0$ and less than four points ramify, then there have to be three, but 
any triply branched cover of ${\bf P}^1$ in characteristic zero has strictly more
than $12(g-1)$ automorphisms (actually, at least $15(g-1)$, as follows easily from the Hurwitz formula).

Since Mumford curves are ordinary (\cite{CKK}, (1.2)), 
the ramification groups at those four points are of the form $\Z_p^{t_i} \rtimes \Z_{n_i}$ for some integers $t_i,n_i$ with $n_i|p^{t_i}-1$ and $i=1,\dots,4$ (this is because 
the second ramification group of an ordinary curve is trivial, cf.\ \cite{Nakajima}). Then the Riemann-Hurwitz formula (including higher ramification groups, cf.\ \cite{Nakajima}) in this case implies 
$$ \sum_{i=1}^4 \frac{1}{n_i}(1-\frac{2}{p^{t_i}}) = -\frac{11}{6}. $$
Suppose that $w$ of the four points are wildly ramified (so $t_i \neq 0$ for them); then the corresponding term on the left hand side is $>1/2$, whereas for the $4-w$ tame points, it is $>-1/2$. Hence the left hand side exceeds 
$w-2$. As it should equal $-11/6$, we get that $w=0$ and all ramification is tame. It is then easy to see that $\{n_i\}=\{2,2,2,3\}$ is the unique solution to the above equation. As the ramification is tame, $p=0$ or $p>3$. \qed

\vspace{1ex}
\paragraph\label{quot} {\bf Trees.} \ The group $\Aut(X)$ is (abstractly) isomorphic to the quotient $N / \Gamma$, where $\Gamma \subseteq PGL(2,k)$ is the Schottky group corresponding to $X$ and $N$ is its normalizer in $PGL(2,k)$. The group $N$ is a finitely generated discrete subgroup of $PGL(2,k)$, and it can be written as the amalgamated product along a certain tree of finite groups (cf.\ \cite{CKK}, \S 2). More precisely, let $\mathcal{T}$ be the Bruhat-Tits tree of $PGL(2,k)$, let ${\mathcal T}_N^*$ be its subtree generated by the limit points and the fixed points of torsion elements of $N$, seen as ends of $\mathcal T$; and let $T^*_N$ be the quotient of ${\mathcal T}_N^*$ by $N$. Then $N$ is the tree product over $T_N^*$, if we make it into a tree of groups by labeling an edge or vertex by the stabilizer of any of its lifts to ${\mathcal T}_N^*$, seen as acted upon by $N$. The ends of $T_N^*$ (which are usually contracted so $N$ becomes a finite tree product) are in bijection with the branch points in $X \rightarrow \Aut(X) \backslash X$, and their stabilizers are the corresponding 
ramification groups.

\vspace{1ex}
\paragraph\label{classical} {\bf Lemma.} \ {\sl Let $X$ be as in \ref{RH}, and $N$ the normalizer of the corresponding Schottky group. Then any finite subgroup of $N$ has order coprime to $p$.}

\vv

\pf It suffices to prove that if we write $N=*_{N_i} N_j$, no group $N_i$ has order divisible by $p$, as any finite subgroup of $N$ is conjugated to some $N_i$. But if this were so, then \cite{CKK} (4.4)(ii) implies that there is wild ramification in $X \rightarrow \Aut(X) \backslash X$, which contradicts \ref{RH}. \qed 

\vspace{1ex}
\paragraph\label{convention} {\bf Convention.} \ We let $D_n$ denote the dihedral group of order $2n$.

\vspace{1ex}
\paragraph\label{amalgams} {\bf Proposition.} \ {\sl Let $X$ be as in \ref{RH}, and $N$ the normalizer of the corresponding Schottky group. Then $N$ is (abstractly) isomorphic to one of the following:
$$ \mbox{(i) } D_2 *_{\Z_2} D_3 \mbox{;\ \ \    (ii) } D_3 *_{\Z_3} A_4 \mbox{;   \ \ \ (iii) } D_4 *_{\Z_4} S_4 \mbox{; \ \ \   (iv) } D_5 *_{\Z_5} A_5, $$
where (iv) only occurs if $p>5$.}

\vv

\pf  We know $p=0$ or $p>3$ by \ref{RH}. To proceed, we use the technique of \cite{CKK}.
Proposition 1 of \cite{CKK} says that if $T^*$ is a subtree of $T^*_N$ having the same ends as $T_N^*$, then $N$ is still the tree product over this $T^*$. Since there are four ends, 
we can assume $T^*$ to be of the form ``two lines (=four ends) connected by a segment''. We know the stabilizers of these ends have to be cyclic groups of order (2,2,2,3) respectively (say, (2,2) occurs  on the left line and (2,3) on the right line). 

We now use the fact that if a finite subgroup $G$ of $N$ stabilizes a vertex, then the edges going out of the vertex have to be stabilized by cyclic groups of order the branching indices of ${\bf P}^1 \rightarrow G \backslash {\bf P}^1$ (at least in the tame case, cf. \cite{CKK} (2.10)). 

A vertex group, say, $G$, is a subgroup of $PGL(2,k)$ of order prime-to-$p$, so by the classification of Dickson (\cite{CKK}, (2.9)), it occurs in the following list:
$$G \in \{ \Z_n, D_n, A_4, S_4, A_5 \}$$ with branching indices of the corresponding ${\bf P}^1$-quotient ${\bf P}^1 \mapsto G \backslash {\bf P}^1$ given by $$(n,n), (2,2,n), (2,3,3), (2,4,4), (2,3,5)$$ respectively. We also know that $p>5$ in the very last case. 

If the segment has an inner vertex, only two edges can emanate from it, hence
by inspection of the above list, its stabilizer should be a cyclic group. Then actually, any stabilizer of a neighbouring vertex should be the same cyclic
group (by the fact that cyclic stabilizers should be maximally cyclic subgroup of corresponding vertex groups, cf. \cite{Herrlich:80}, Lemma 1 and \cite{CKK} (4.1)). Now the next edge has to be stabilized by the same $\Z_n$ (corresponding to the second ramification point in the $\Z_n$-quotient of ${\bf P}^1$). Proceeding in the same way, all stabilizers of edges and vertices on the inner segment are the same cyclic group, so we do not change the
abstract structure of the corresponding tree product by replacing the segment
by just one edge with this cyclic group. Thus, $T_N^*$ looks as in the following picture (where arrows indicate ends): 

\vv

{\footnotesize
\unitlength0.7pt
\begin{center}
\begin{picture}(40,30)
\multiput(0,15)(35,0){2}{\circle*{4}}
\multiput(0,15)(35,0){1}{\line(1,0){35}}
\multiput(0,15)(35,0){2}{\vector(0,1){20}}
\multiput(0,15)(35,0){2}{\vector(0,-1){20}}
\put(-3,37){$2$}
\put(-3,-15){$2$}
\put(32,37){$2$}
\put(32,-15){$3$}
\put(-15,12){$U$}
\put(40,12){$V\ \ \ \ .$}
\put(12,3){$H$}
\end{picture}
\end{center}}

\vv

The possibility that $U$ is cyclic can be excluded, because maximality of cyclic subgroups then implies that $H$ is trivial, but then there is no choice for $V$ that leads to the desired stabilizers of ends (2,3).  

Further inspecting the above list, we see that only $U=D_n$ is possible (so the left line is stabilized by elements of order two). Then the middle segment has to be stabilized by $H=\Z_n$. If the right vertex also has a dihedral group, then it has to be $V=D_3$ (because we need a (2,3)-ramification) and $n=2$, so we are in case (i). 
Similarly, if the right vertex is $A_4, S_4$ or $A_5$, then $n=3,4,5$ and we are in case (ii), (iii) and (iv), respectively. \qed

\vv

\vv

\sectioning\label{lame}\ 
{\bf Classification of groups of Lam\'e type}

\vv

\vspace{1ex}
\paragraph\label{lame1} {\bf Reminder.} \ Let $X$ be a Mumford curve of genus $g \in \{5,6,7,8\}$ with automorphism group $G:=\Aut(X)$ of order $12(g-1)$
(which is $48, 60, 72$ and $84$, respectively). We now know that $G$ is a quotient of some group $N$ as in proposition \ref{amalgams} by a free group $\Gamma$ of rank $g$. Our aim in this paragraph is to systematically determine
all possible such $G$ from these facts. We start with an obvious lemma:

\vspace{1ex}
\paragraph\label{converse} {\bf Lemma.} \ {\sl Let $N$ be a finitely generated discrete subgroup of $PGL(2,k)$ that is abstractly isomorphic to an amalgam $U*_Z V$ of finite groups, and assume that we are given a surjective group homomorphism $\phi \ : \ N \rightarrow G$ to some finite group $G$. Then if the kernel $\Gamma:=\ker(\phi)$ is free of finite type, it maps $U$ and $V$ isomorphically into $G$, and their images generate $G$.} 

\vv

\pf The morphism $\phi$ restricted to $U$ is injective, since an element of $U$ in the kernel would give an element of finite order in $\Gamma$, which contradicts the fact that $\Gamma$ is a free group. Since any element of $N$ is already a word in
letters from $U$ and $V$ and $\phi$ is surjective, any element of $G$ is a
word in the images of letters from $U$ and $V$. \qed

\vv

In the next few lemmas, we dismiss most possibilities for $G$ using the elementary theory of groups, applying lemma \ref{converse}. For the reader to appreciate the drastic result that only two groups remain, let us remark that there are 52, 13, 50 and 15 non-isomorphic groups of order 48, 60, 72 and 84, respectively. 

Most of the arguments involve counting Sylow groups and a discussion of split extensions. As a compromise to readibility, we have dismissed the two most technical such arguments to an appendix, and present here only the shorter ones.

Some of the claims below can also be checked by computer, and the authors have used the system {\tt GAP} (\cite{GAP}) for double checking a few of the results.

\vspace{1ex}
\paragraph\label{case1} {\bf Lemma.} \ {\sl In the situation of the main theorem, if $N \cong D_5 *_{\Z_5} A_5$, then $g=6, p\neq 5$ and $G \cong A_5$.}

\vv

\pf By the previous lemma, $G$ has a subgroup isomorphic to $A_5$, hence the order of $G$ has to be 60, and $g=6$. Also, $p \neq 5$ by Dickson's classification. \qed

\vspace{1ex}
\paragraph\label{case2} {\bf Lemma.} \ {\sl In the situation of the main theorem, if $N \cong D_4 *_{\Z_4} S_4$, then $g=5$ and $G \cong S_4 \times \Z_2$.}

\vv

\pf Let $U$ and $V$ denote the respective isomorphic images of $D_4$ and $S_4$ in $N$. As $G$ contains $V$, its order has to be 48 or 72, and, correspondingly, $g=5$ or $g=7$. 

\vv

(\ref{case2}.1) {\sl The case $g=5$.} \ $V$ is of index two in $G$, hence normal, and as $U \not \subset V$ (since otherwise, $G=V$) and $U \cap V$ contains all elements of order four in $U$, there is an element of exact order two outside $V$, so the sequence $1 \rightarrow V \rightarrow G \rightarrow \Z_2 \rightarrow 1$ splits. Since $\Aut(S_4)=S_4$ and $S_4$ has two conjugacy classes of elements of order two, there are at most two such non-isomorphic non-trivial split extensions, and if we write $G=S_4 \rtimes \Z_2$, $\Z_2 = \langle \gamma \rangle$, they are given by  $\gamma \mapsto (12)$ and $\gamma \mapsto (12)(34)$. However, both of these are isomorphic to the
trivial one by multiplying on the left with $(12)$ and $(12)(34)$, respectively.
Hence $G = S_4 \times \Z_2$, and then there does exist a surjective homomorphism $N \mapsto G$ (the image of $D_4$ is $\langle (12)(34)\gamma, (1234) \rangle$, and $Z$ maps to $\langle (1234) \rangle $). \

\vv

(\ref{case2}.2) {\sl The case $g=7$.} \ The Sylow theorem implies that there are either one or four 3-Sylow groups. Suppose there is only one, say, $P$ (of order 9). then $G$ contains at most four subgroups of order 3 (as all are contained in $P$), but already $V$ contains four such groups, hence $P \subseteq V$, so $9|24$, a contradiction. 

Hence $G$ has precisely four 3-Sylows, say, $\{P_i\}_{i=1}^4$. $G$ acts on these by conjugation, so gives a homomorphism $\phi: G \rightarrow S_4$, and $\ker(\phi)=\bigcap_{i=1}^4 N(P_i)$, where $N(-)$ denotes normalizer. This kernel is non-trivial, since $|G|/|S_4|=3$. For all $i$, $N(P_i)$ contains 18 elements. Since the $N(P_i)$ cannot have any $P_i$ (of order nine) in common, $\ker(\phi)$ contains three or six elements. The restriction of $\phi$ to $V$ gives a homomorphism $V \cong S_4 \rightarrow A_4$ whose kernel has order dividing 6. But $V \cong S_4$ does not have normal subgroup of order 2, 3 or 6, so $V \cap \ker(\phi)=\{1\}$ and hence $|\ker(\phi)|=3$.

Therefore, the exact sequence $1 \rightarrow \Z_3 \rightarrow G \rightarrow S_4 \rightarrow 1$ splits. The semi-direct product is non-trivial, since otherwise, there would be no copies of $D_4$ outside $S_4$ in $G$. There is a unique such non-trivial semi-direct product given by the signature homomorphism $S_4 \mapsto \Aut(\Z_3)=\Z_2$. If $\gamma \in \Z_3$ and $\alpha \in S_4$, we have $\gamma \alpha = \alpha \gamma^{\mbox{{\tiny sgn}}(\alpha)}$ in $G=\Z_3 \rtimes S_4$. 

The number of 2-Sylows of $G$ is 1, 3 or 9, but there are already three in $V$ and one more equals $U$, so there must be 9. For every two-Sylow $Q$ of $V \cong S_4$ (of which there are three), we get three two-Sylows of $G$ (namely, the conjugates by $\langle \gamma \rangle$) and these are all two-Sylows of $G$. Now all of these nine 2-Sylows of $G$ are isomorphic to $D_4$, and hence these are all the subgroups of $G$ isomorphic to $D_4$.  One computes immediately that the 2-Sylows of $S_4$ intersect in a four group, and then the same holds for all $D_4$ in $G$. Hence no $D_4$ in $G$ can share a cyclic $\Z_4$ with $V$, and so there is no map $N \rightarrow G$. \qed

\vspace{1ex}
\paragraph\label{case3} {\bf Lemma.} \ {\sl In the situation of the main theorem, if $N \cong D_3 *_{\Z_3} A_4$, then $g=6$ and $G=A_5$.}

\vv

\pf This time, all $g$ can occur. Let $U$ and $V$ denote the respective isomorphic images of $A_4$ and $D_3$ in $N$.

\vv

(\ref{case3}.1) {\sl The case $g=5$.} \ This case is dismissed by the theorem \ref{thm-groupthm} in the appendix.

\vv

(\ref{case3}.2) {\sl The case $g=6$.} \ If there were a unique  5-Sylow $P$ in $G$, then $G/P \cong A_4$ and as $P \cap D_3 = \{ 1\}$, we would get $D_3 \subseteq A_4$ in the image, a contradiction. 
Therefore, $G \cong A_5$, as is shown in the appendix \ref{thm-groupthm2}.

\vv

(\ref{case3}.3) {\sl The case $g=7$.} \ See appendix \ref{thm-groupthm}.

\vv

(\ref{case3}.4) {\sl The case $g=8$.} \ Then $G$ has a unique normal 7-Sylow group $P$, and $U \cong A_4$ and $V \cong D_3$ should still be subgroups of $G/P$, hence $G/P \cong A_4$, and so $D_3 \subseteq A_4$, a contradiction. \qed 

\vspace{1ex}
\paragraph\label{case4} {\bf Lemma.} \ {\sl In the situation of the main theorem, if $N \cong D_2 *_{\Z_2} D_3$, then $g=5$ and $G=S_4 \times \Z_2$.}

\vv

\pf Again, all $g$ can occur. Let $U$ and $V$ denote the respective isomorphic images of $D_2$ and $D_3$ in $N$. 
\vv

(\ref{case4}.1) {\sl The case $g=5$.} \ See appendix \ref{thm-groupthm}.

\vv

(\ref{case4}.2) {\sl The case $g=6$.} \ If there is a unique 5-Sylow  $P$, then the quotient is a group of order 12 in which $D_3$ is normal (since of index two); hence it has a unique 3-Sylow, and pulling it back to $G$ we get a normal subgroup of order 15, hence a $\Z_{15}$, hence a unique 3-Sylow in $G$, which intersects $D_3$ in $\Z_3$. The quotient of $G$ by this 3-group should have order 20 and be generated by $D_2$, a contradiction.
Hence the number of 5-Sylows of $G$ is 6. The appendix \ref{thm-groupthm2} shows that  $G \cong A_5$.

\vv

(\ref{case4}.3) {\sl The case $g=7$.} \ See appendix \ref{thm-groupthm}.

\vv

(\ref{case4}.4) {\sl The case $g=8$.} \ There is a unique normal 7-Sylow group $P$; let $\pi : G \rightarrow G/P$ be the corresponding quotient map. 

The number of 3-Sylows is 1, 3 or 7. If there are 7, their normalizers are
of order 12, so map by $\pi$ injectively into $G/P$ (which is of order 12). hence $\pi$ gives a splitting $G \cong \Z_7 \rtimes M$ for some such normalizer $M$. The product is not direct as $G$ cannot have a quotient $\Z_7$ (it is generated by $D_2$ and $D_3$). The extension corresponds to a map $M \rightarrow \Aut(\Z_7)=\Z_6$. If $M=D_6$, there is only one such non-trivial homomorphism with image $\Z_2$, and then $G$ contains a normal $\Z_7 \times \Z_6$. So there
is a quotient $G/\Z_6$ of order 14, which cannot be generated by $D_2$ and $D_3$. The second possibility is $M \cong \Z_3 \rtimes \Z_4$, and this leads
similarly to a normal $\Z_7 \times \Z_6$. 

If there are 4 3-Sylows, then the kernel of the action $G \rightarrow S_4$ has 7 or 21 elements, but if there were 21, $G$ would have one or seven 3-Sylows. Hence the kernel has 7 elements, so equals $P$, and the image of this map is $A_4$, which cannot contain $D_3$. 

Hence there is a unique 3-Sylow, which intersects $D_3$ in a $\Z_3$, so the quotient of $G$ by this 3-Sylow is a group of order 28 generated by $D_2$, a contradiction. \qed

\vv

These lemmas finish the proof of the main theorem. \qed

\vv

\vv

\sectioning \label{realization} \ {\bf Realization as automorphism groups.}

\vv

\vspace{1ex}
\paragraph\label{introd} {\bf Introduction.} \ In this section, we show that the possible pairs $(N,G)$ that are left by the group theoretical arguments of the previous paragraph can indeed be realized
by Mumford curves. We note that it suffices to prove that such $N$ 
can be realized as discrete subgroups of $PGL(2,k)$, and for this, one 
could use known methods (e.g., Herrlich, \cite{Herrlich:80b}). However, we will take a different approach as follows. 
For a given $g$ and $G$, we first explicitly construct the family of algebraic  curves
of genus $g$ whose automorphism group contains $G$. This is not too difficult as $g$
is not too large. We then calculate stable models for the bad fibers of 
these families, and this allows us to read off the Mumford loci and the 
corresponding normalizers of Schottky groups. Recall that we only have to do this for the remaining case of proposition A.

\vspace{1ex}
\paragraph\label{prop} {\sl Proof of Proposition A.} \ Let $X$ be an algebraic curve with $g=5$ and $G=\Aut(X)=S_4 \times \Z_2$. Let $Y$ be the quotient of $X$ by $\Z_2$, which is an $S_4$-Galois cover of ${\bf P}^1$, whose ramification is tame above at most four points with indices taken from $(2,2,2,3)$. One checks immediately using the Riemann-Hurwitz formula for $Y \mapsto {\bf P}^1$ that the only possibility is that it ramifies over four points with indices $(2,2,2,3)$, and then, $Y$ is of genus three. This implies that $X$ is an unramified cover of $Y$. 

Now all curves of genus three with an octahedral automorphism group (or larger) have been determined by Wiman in \cite{Wiman}, using the fact that the canonical embedding maps such curves to the plane with a number of natural sets of fixed points (inflection, double tangent and sextactic). He finds exactly the one-parameter family $C^\sim_\alpha$ as in proposition A, and observes that 
the bitangents $x \pm y \pm z = 0$ touch $C^\sim_\alpha$ in 8 points. One calculates immediately that these fall naturally into two $S_4$-orbits, one of which is given by the four points $$ \{p_i\}_{i=1}^4 = \{ (1:\omega:\omega^2), (1:\omega:-\omega^2),(1:-\omega:-\omega^2),(1:-\omega:\omega^2)\},$$
where $\omega$ is a fixed third root of unity. The other $S_4$-orbit is given by $\{p'_i\}=\{p_i^\sigma\}$ for $\langle \sigma \rangle = \mbox{Gal}({\bf Q}(\omega)/{\bf Q})$. 

Since $X$ is an unramified 2-cover of $Y=C^\sim_\alpha$, it should correspond to an $S_4$-invariant two-torsion point in Jac$(Y)$, i.e., a degree zero $S_4$-invariant divisor on $Y$. Now the two lines $L_1,L_2=0$ through $(p_1,p_2)$ and $(p_3,p_4)$ together with the two lines $L_3,L_4=0$ through $(p_1,p_3)$ and $(p_2,p_4)$
generate a pencil of conics $P_\lambda \ : \ L_1 \cdot L_2 + \lambda \ L_3 \cdot L_4 = 0$ through all four $\{p_i\}$, and it is tangent to one (hence, all) of the bitangents above precisely if $\lambda=\omega^2$. Taking Galois conjugates gives a similar pencil $P'_\lambda$, and so finally 
$ P''_\lambda : P_\omega + \lambda P'_{\omega^2} = 0$ gives a pencil of conics 
that are tangent to the four lines $x \pm y \pm z = 0$. 

Now suppose $\alpha \neq -1$. Then $C^\sim_\alpha$ does not intersect $P''_\lambda$ at its four base points $(1:\pm 1:\pm 1)$, so there is a well defined 
function $f \ : \ C^\sim_\alpha \rightarrow {\bf P}^1: P \mapsto \lambda_P$, where $\lambda_P$ is
the value of $\lambda$ such that $P''_\lambda$ goes through $P$. This function $f$ has divisor div$(f)=2(\sum_{i=1}^4 (p_i-p_i'))$, as $\lambda=0$ (respectively, $\lambda=\infty$) give precisely the conics tangent at $\{ p_i \}$ (respectively, $\{p'_i\}$). Finally, $D=\sum_{i=1}^4 (p_i-p_i')$ is a two-torsion element in Jac$(C^\sim_\alpha)$, which is clearly $S_4$-invariant and non-trivial (since the four points $p_i$ are not collinear). Hence we have found our desired family $C_\alpha$.

Note that if $\alpha=-1$, we a priori have such a map $f$ defined outside
the base points of the pencil $P''_\lambda$, but then it can be extended 
to all of $C^\sim_{-1}$.  

We still have to prove that we have found all such curves, i.e., that there
is no further $S_4$-invariant 2-torsion point in Jac$(C^\sim_\alpha)$. However, 
this will be shown at the end of the proof using one of the degenerate fibers
which we now study first. 

\vv

Those degenerate fibers of the family $C^\sim_\alpha$ occur precisely at $\alpha\in\{\pm2, \infty\}$ (easily seen calculating singularities) --- as a matter of fact, the singular fibers of $C_\alpha$ occur at precisely the same values of $\alpha$, since the cover is \'etale. We calculate as in \cite{Kato:01}, \S 3.
It is known that $C^\sim_\alpha$ is a Mumford curve precisely when $\alpha$ 
coincides with a degeneration point over the residue field $\bar{k}$, such that the fiber over that degeneration point has a stable model which is the union of rational curves intersecting in rational point. We will prove that this is the case for $\alpha=\infty$ and $\alpha=-2$, whereas $\alpha=2$ has (potentially) 
good reduction. 

\vv

(\ref{prop}.1) \ If $\alpha=\infty$, then the special fiber is a rational curve with three nodes, which, after blowing up the double points, becomes a line (solid in the picture) intersecting each of three lines (the special fibers; dashed in the picture) in exactly two points. 

\begin{center}\scalebox{.6}{\includegraphics{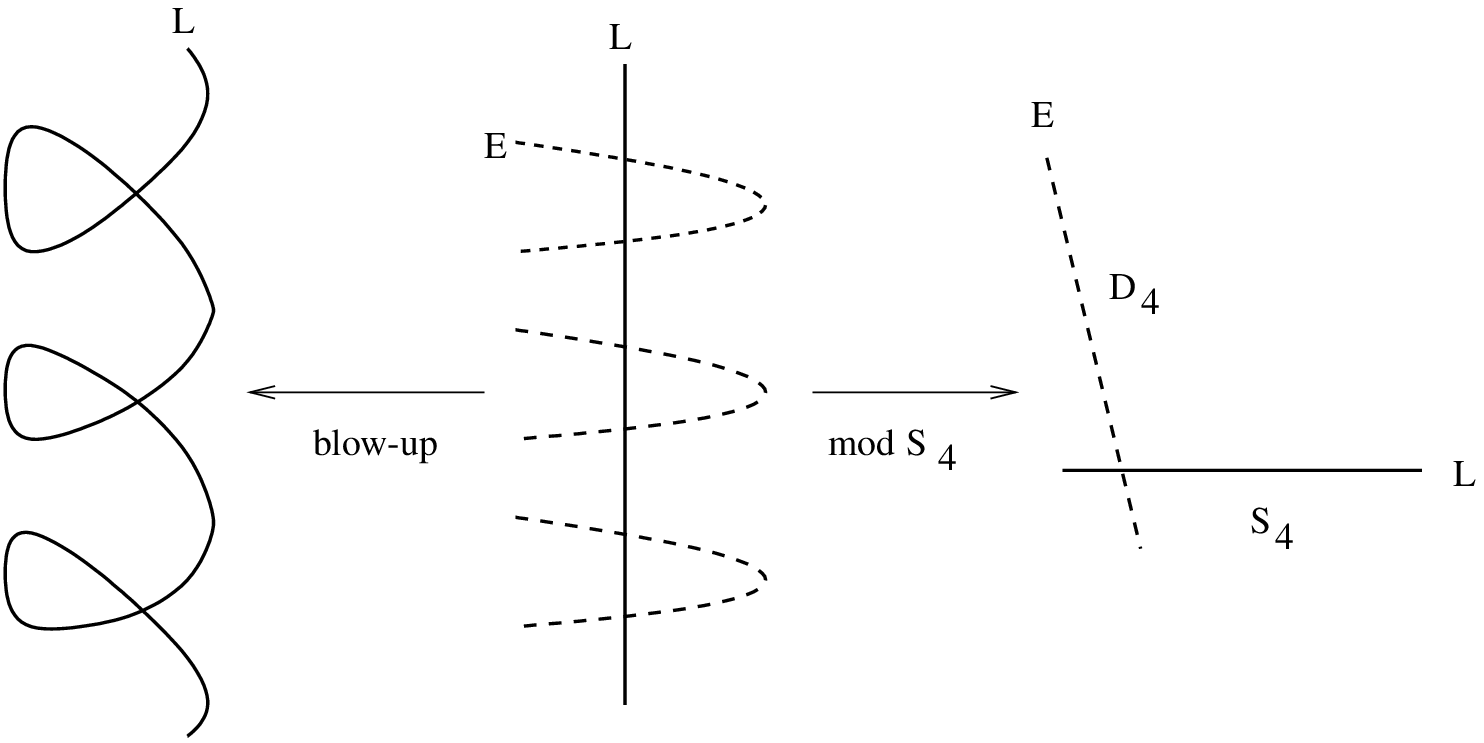}}
\end{center}

Now $S_4$ acts without inversion on the intersection dual graph of this configuration. The original line is 
stabilized by $S_4$, and ones sees that the exceptional divisors are stabilized by $D_4$ (the only group of index three in $S_4$ up to isomorphism), and these two groups intersect in a $\Z_4$. Hence if we quotient by $S_4$, we get two lines (the image of the original line and of the exceptional divisor) which are stabilized by $S_4$ and $D_4$, respectively, and so that the intersection point is stabilized by $\Z_4$.
So the quotient of a stable model of $C^\sim_\alpha$ over the valuation ring of $k$ for $|\alpha|>1$ is a rational curve with intersection graph of the special fiber equal to the tree product $S_4 *_{\Z_4} D_4$, and so we are in case (a1). 

Now we pass to $C_\alpha$. The support of the divisor $D$ restricts to smooth points of the rational curve $L$. Therefore, a semi-stable model of $C_{\infty}$ is given by two rational curves connected by 6 non-intersecting exceptional 
divisors (which are rational lines) --- the intersection dual of this consists of two points connected by six edges. Thus, one calculates in the same way 
that case (a1) occurs.

\vv

(\ref{prop}.2) \ If $\alpha=-2$, the degenerate fiber consists of 4 lines intersecting pairwise, each of which is stabilized by a $D_3$ (index 4 in $S_4$), and on which $S_4$ acts transitively. If we blow up all 6 double points, $S_4$ acts without inversion on the intersection graph, and the quotient consists of two lines with one double point (the image of a line and the image of an exceptional divisor). 

\begin{center}\scalebox{.6}{\includegraphics{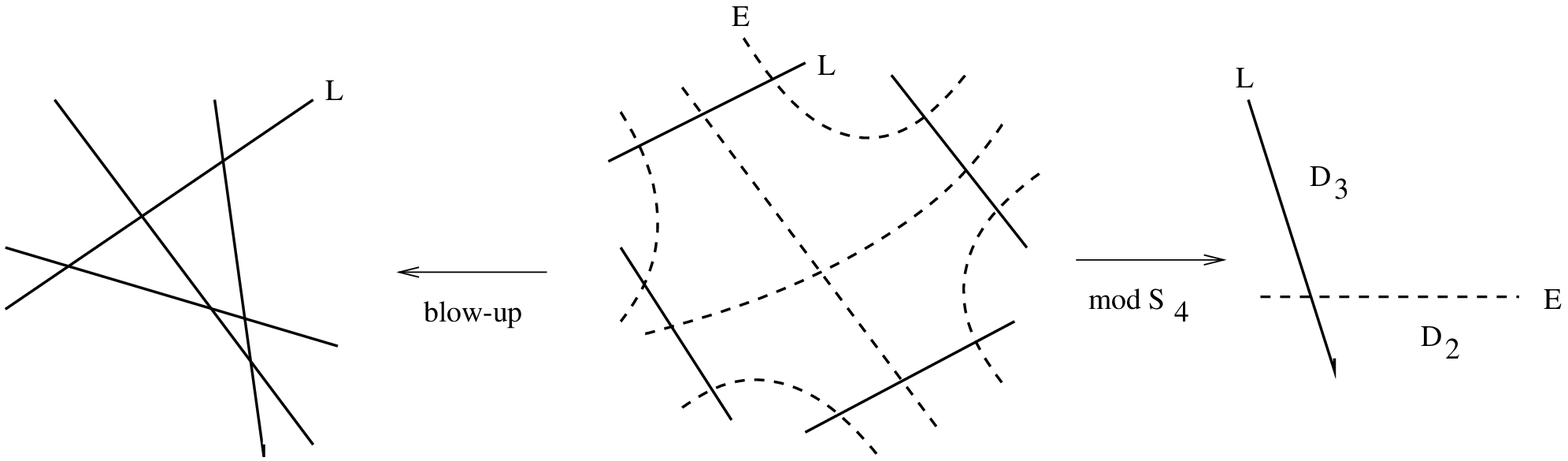}}
\end{center}

The image of a line is stabilized by $D_3$ and the exceptional divisor by $D_2$ (the stabilizer of a double point before the blow-up -- note that it is a group of index 6 in $S_4$ all of whose elements have to have order two, since they should interchange the two lines through the double point), and these intersect 
in a $\Z_2$. Hence the quotient of a stable model of $C^\sim_\alpha$ over the valuation ring of $k$ for $|\alpha+2|<1$ is a rational curve with intersection graph of the special fiber equal to the tree product $D_3 *_{\Z_2} D_2$, and so we are in case (a2). 

Again, the support of the divisor $D$ restricts to smooth points, and so passing from $C^\sim_{-2}$ to $C_{-2}$ just causes the solid lines in the above picture to be doubled. Again, we see that we are in case (a2).

\vv

(\ref{prop}.3) \ We will now show that curves around $\alpha=2$ have good reduction. In our model
of the pencil $C^\sim_\alpha$, we see only a double conic at $\alpha=2$, so we change
the pencil by picking two new generators to $$ C^*_t \ : \ (x^2+y^2+z^2)^2 + t (x+y+z)(x-y+z)(x+y-z)(x-y-z) = 0$$ 

We first look at what happens \'etale locally around a bitangent point (in the dashed circle below), where
the conic $x^2+y^2+z^2=0$ and the line $x+y-z=0$ can be replaced by $X=0$ and $Y=0$, and the pencil looks like $X^2-tY=0$. 

\begin{center}\scalebox{.6}{\includegraphics{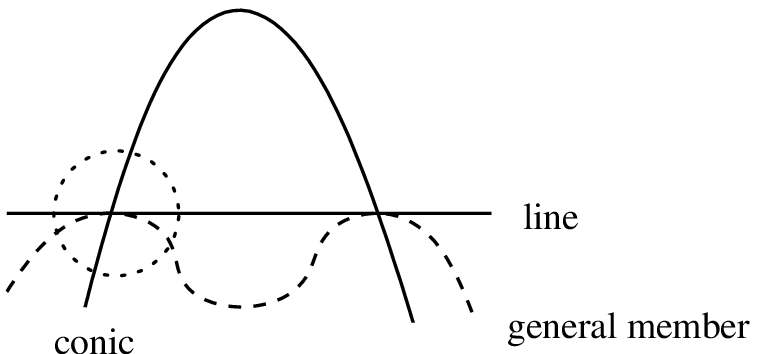}}
\end{center}

To split the multiple fiber, set $t=\tau^2$, then the corresponding surface
in $(X,Y,\tau)$-space $X^2-\tau^2 Y=0$ has a one-dimensional singularity along $X=\tau=0$, which blows up to an exceptional divisor $\tau^2=0$ and a conic $(\frac{X}{\tau})^2-Y=0$. 

If we are not around a bitangent point, setting $t=\tau^2$ just leads to two
transversally intersecting components.

All of this shows that, in the pencil $C^*_t$, blowing up along 
$x^2+y^2+z^2$ gives a curve which is a double cover of this conic ramified 
over eight (bitangent) points, so is a hyperelliptic curve of genus three, 
so we have good reduction.  

\vv

Finally, we turn to the proof that the divisor $D$ is the unique $S_4$-invariant two torsion point in $J_\alpha:=\mbox{Jac}(C^\sim_\alpha)$. For this, notice that 
$J_\alpha[2]$ is a locally {\sl constant} group scheme $(\underline{\Z/2\Z})^6$, hence to show
that there is a unique $S_4$-invariant vector in this space, we can check
it at any particular fiber. We choose $\alpha=2$, then by the calculation in the previous paragraph the curve
is a hyperelliptic cover of ${\bf P}^1$ ramified above precisely eight points,
corresponding to the eight bitangent points $\{ p_i, p'_i \}_{i=1}^4$, and this precribes the action of $S_4$ on them. Now the 2-torsion in the Jacobian
of a hyperelliptic curve like $C^\sim_2$ is well know to be spanned by divisors
supported at those ramification points (=Weierstrass points). Suppose $D'$ is
an $S_4$-invariant element in $J_2[2]$. We know supp$(D') \subseteq \{ p_i, p'_i \}_{i=1}^4$. If a $p_i$ occurs, say, $p_1$, then its whole $S_4$-orbit $\{ p_i \}$ should occur with the same multiplicity. As the degree of $D'$ should be zero, an element outside this orbit should also occur, hence some $p_i'$, hence
its whole $S_4$-orbit, so $D' = n \sum_{i=1}^4 (p_i-p'_i)$ for some odd $n$, which is the same element of $J_2[2]$ as $D$. This finishes the proof of the proposition.
\qed

\vv

\vv

\renewcommand{\thesection}{\Alph{section}}
\setcounter{section}{0}

{\footnotesize \sectioning\label{appendix}\  {\bf Appendix: more group theory}

\vv

\paragraph\label{thm-groupthm} {\bf Theorem.} \
{\sl Let $G$ be a finite group containing two subgroups $U$ and $V$ which generate
$G$. Consider the following cases:

}(a) {\sl $U\cong D_3$, $V\cong D_2$, and $U\bigcap V\cong\Z/2\Z$, 

}(b) {\sl $U\cong D_3$, $V\cong A_4$, and $U\bigcap V\cong\Z/3\Z$.

Then

$(1)$ If $|G|=48$ and (a) holds, $G$ is isomorphic to $S_4\times\Z/2\Z$, whereas \textup{(b)} cannot hold,

$(2)$ there is no group $G$ with $|G|=72$ such that either \textup{(a)} or \textup{(b)} holds.}

\vv

\pf We will use the following presentation of the groups: set  $G=\langle x,y,z\rangle$ and

in case (a): $U=\langle x,y\,|\,x^3=y^2=1,\ yxyx=1\rangle$, $V=\langle y,z\,|\,
y^2=z^2=1,yz=zy\rangle$, $U\bigcap V=\langle y\rangle$;

in case (b): $U=\langle x,y\,|\,x^3=y^2=1,\ yxyx=1\rangle$, $V=\langle x,z\,|\,
x^3=z^2=(zx)^3=1\rangle$, $U\bigcap V=\langle x\rangle$.

\vspace{1ex}
\paragraph\label{lem-1} {\bf Lemma.} \ 
{\sl The $2$-Sylow subgroups of $G$ are not normal.}

\vv

\pf
Let $P_2$ be a $2$-Sylow subgroup of $G$.
The elements $y$ and $xy$ are of order $2$, while 
the product $xy\cdot y=x$ is of order $3$. If $P_2$ were normal, hence
unique, both $y$ and $xy$ would belong to $P_2$, and hence $x\in P_2$, which
is absurd.
\qed

\vspace{1ex} \paragraph \label{lem-2} {\bf Lemma.} \
{\sl The $3$-Sylow subgroups of $G$ are not normal.}

\vv

\pf
Suppose a $3$-Sylow subgroup $P_3$ is normal. 

Case (a): Since $x\in P_3$, we have $G/P_3=\langle yP_3,zP_3\rangle
\cong D_2$, and hence $|G/P_3|=4$. It follows that the order of $G$ must be
of the form $2^2\cdot 3^k$, which is not the case when $|G|=48$ or $72$.

Case (b): In $A_4$ there are two elements of order $3$ of which the 
product is of order $2$; e.g., $(123)(124)=(13)(24)$. This implies $P_3$ contains
an element of order $2$, which is a contradiction. \qed

\vspace{1ex} \paragraph \label{lem-3} {\bf Lemma.} \ {\sl
There are precisely four $3$-Sylow subgroups in $G$.}

\vv

\pf
Case $|G|=72$: Immediate from Lemma \ref{lem-2} and Sylow's theorem.

{ Case $|G|=48$:} The number of $3$-Sylow
subgroups is either $1$, $4$, or $16$. We know already it is not $1$. Suppose 
it is $16$. Then there are $2\cdot 16=32$ elements of order $3$, and hence all 
the rest have to form the unique $2$-Sylow subgroup, violating Lemma \ref{lem-1}.
\qed

\vv

To proceed, we fix the transitive representation $\rho\colon G\rightarrow S_4$,
induced from the action of $G$ on the set of all $3$-Sylow subgroups.
We denote its kernel by $K$.

\vspace{1ex} \paragraph \label{lem-4} {\bf Lemma.} {\sl
The subgroup $K$ does not contain a $3$-Sylow subgroup. Moreover, we have
$|\rho(G)|=12$ or $24$.}

\vv

\pf
Suppose a $3$-Sylow subgroup $P_3$ is contained in $K$. Then elements of
$P_3$ normalize another $3$-Sylow subgroup $P'$,
$P_3\neq P'$. This means that $G$ contains a subgroup isomorphic to a semi-direct
product of $P_3$ and $P'$, which is a $3$-group of order strictly larger than 
that of $P_3$, whence a contradiction. Hence the first assertion is proved.
It then follows that $|\rho(G)|$ is divisible by $3$ since, otherwise, $\rho$
would map each $3$-Sylow subgroup to the unit element.
Hence $\rho(G)$ contains an element of order $3$, which implies in particular that
$\rho(G)$ is doubly transitive. Now $|\rho(G)|$ is $4\cdot 3$ times the order of
the stabilizer of two points, whence the result. 
\qed

\vspace{1ex} \paragraph \label{lem-5} {\bf Lemma.} \ {\sl
$x\not\in K$.}

\vv

\pf
If $|G|=48$, then $|K|=2$ or $4$, hence the assertion.
Suppose $|G|=72$ and $x\in K$. 
Since in this case $|K|=3$ or $6$, $\langle x\rangle$ is the unique $3$-Sylow 
subgroup of $K$, hence is normal in $G$.
In Case (b), this contradicts the fact that in $V$, $3$-Sylow subgroups are not
normal. In Case (a), we would have $G/K=\langle yK,zK\rangle\cong D_2$, which
contradicts Lemma \ref{lem-4}.
\qed

\vspace{1ex}\paragraph \label{lem-6} {\bf Lemma.} \ {\sl
The order of $\rho(G)$ is actually $24$. In particular, $|K|=2$ (if $|G|=48$)
or $3$ (if $|G|=72$).}

\vv

\pf
Suppose $|\rho(G)|=12$, i.e., $\rho(G)\cong A_4$.
Since $x\not\in K$ (Lemma \ref{lem-5}), we may assume $\rho(x)=(123)$.
It holds that $yxy=x^{-1}$, but there is no element $\sigma$ of order $2$ in $A_4$ 
such that $\sigma(123)\sigma=(132)$. This is a contradiction.
\qed

\vv

{\sl Proof of Theorem \ref{thm-groupthm}}. \ 
Since $G=\langle x,y,z\rangle$, we know $\rho(G)=\langle\rho(x),\rho(y),
\rho(z)\rangle=S_4$.
As in the proof of Lemma \ref{lem-6}, we may set $\rho(x)=(123)$, and, since 
$\rho(y)$ has to normalize $(123)$, $\rho(y)=(12)$.
Since $\rho(G)=S_4$, $\rho(z)$ has to be of order $2$, which commutes with 
$(12)$. Such an element in $S_4$ is either $(12)$, $(34)$, or $(12)(34)$, where
the first one can be dismissed since $\rho(G)=S_4$. We can without loss of generality suppose
$\rho(z)=(12)(34)$, since if $\rho(z)=(34)$, we can replace $z$ by $yz$.
Set $W=\langle x,z\rangle$.
Note that $\rho(W)=\langle(123),(12)(34)\rangle=A_4$, hence in particular, 
$y\not\in W$ (otherwise, the image of $\rho$ would be $A_4$).
We set $K=\langle w\rangle$, where $w$ is of order $2$ or $3$ for $|G|=48$ or $|G|=72$ respectively.

\medskip
{\sl Case} (a) {\sl with} $|G|=48$ ($|K|=2$): 
Since $K$ is normal in $G$ and of order $2$, we have $K\subseteq Z(G)$.
We claim that $K$ is contained in $W$.
Indeed, since $y$ normalizes $W$ (recall: $yxy=x^{-1}$ and $yzy=z$), $G$ is a
semi-direct product of $W$ and $\langle y\rangle$, and hence $|W|=24$.
Since $|\rho(W)|=12$, we have $K\subset W$.

Now set $u=x^{-1}zx$ and $v=x^{-1}ux$. We have
$\rho(u)=(13)(24)$ and $\rho(v)=(14)(23)$.
In particular, $\rho(zuv)=1$, whence $zuv\in K=\langle w\rangle$.
If $zuv=1$, then $\{1,z,u,v\}$ forms a subgroup normalized by $x$, which would imply
that $W$ is a semi-direct product of $\langle x\rangle$ and
$\{1,z,u,v\}$, hence of order $12$, which is a contradiction.
Hence $zuv=w$. Set $r=zw$, $s=x^{-1}rx$, $t=x^{-1}sx$. Then we have
$rst=1$.
Hence we obtain a subgroup $\{1,r,s,t\}$ of order four normalized by $x$. 
Hence $\langle x,r\rangle\cong A_4$.

Now since $w\in K\subseteq Z(G)$, one sees easily that $y$ normalizes $\langle
x,r\rangle$. Hence $\langle x,y,r\rangle$ is of order $24$. 
Since $\rho(r)=\rho(z)$, and since $\langle\rho(x),\rho(y),\rho(z)\rangle=S_4$, 
$\rho$ gives the isomorphism $\langle x,y,r\rangle\cong S_4$. 
Since $w\not\in\langle x,y,r\rangle$, and $w\in Z(G)$, we have
$G\cong S_4\times\langle w\rangle\cong S_4\times\Z_2$.

\medskip
{\sl Case} (a) {\sl with} $|G|=72$ ($|K|=3$):
An easy calculation shows $\rho((xz)^3)=1$, which implies $(xz)^3=1$, $w$, or
$w^{-1}$.
If $(xz)^3=1$, then we see $W\cong A_4$, and, since $y$ normalizes
it, we would have $|G|=|\langle x,y,z\rangle|=24$, which is absurd.
Replacing $w$ by its inverse if necessary, we may suppose $(xz)^3=w$.
Since $\Aut(\Z_3)\cong\Z_2$, we have $xwx^{-1}=w$.
Since $K$ is normal, $zwz=w$, or $w^{-1}$.
We claim $zwz=w$. Indeed, otherwise, $w=(xz)^{-3}w(xz)^3=w^{-1}$, contradiction.
Hence $K$ is contained in the center of $H=\langle x,z,w\rangle$.
Note that $\rho(H)=\rho(W)=A_4$.

Let $N$ be the unique normal subgroup of order 4 (and hence index 3) in $A_4$, and set $L=\rho^{-1}(N)$.
Then $L$ is a subgroup of $H$ of index $3$, hence of order $12$.
Let $S$ be a $2$-Sylow subgroup of $L$. Then $S$ is isomorphic to $N\cong
(\Z_2)^2$.
Since $K$ is included in the center of $H$, it follows that $L\cong K\times S$.
Since $L$ is a normal subgroup of $H$, and $S$ is the unique $2$-Sylow subgroup
of $L$, we deduce that $S$ is normal in $H$, and hence, $S$ is the unique 
$2$-Sylow subgroup of $H$. In particular, $z\in S$.

Next consider $\langle x,S\rangle\subset H$.
In this, $S$ is a normal subgroup, and hence $\langle x,S\rangle$ is a semi-direct
product of $\langle x\rangle$ and $S$, so it is of order $12$.
Now, since $z\in S$, we see $W\subseteq\langle x,S\rangle$, 
while $\rho(W)=A_4$ is also of order $12$.
Hence $W=\langle x,S\rangle$, so
$|W|=12$, which would imply that, since $y$ normalizes
$W$, $|\langle x,y,z\rangle|=24$, which is a contradiction.

\medskip
{\sl Case} (b) {\sl with} $|G|=48$ ($|K|=2$):
Since $\rho(yzyz)=1$, $yzyz=1$ or $w$.
If $yzyz=1$, then the group $G$ with $U$ and $V'=\langle y,z\rangle$ fits in 
with Case (a), hence is isomorphic to $S_4\times\Z_2$. But then, there is
no subgroup $V$ of $G$ as in the assumption. Hence we have $yzyz=w$.
Since $K\subseteq Z(G)$, we have $yzy=zw$.
Then, $1=y(xz)^3y=(yxyyzy)^3=(x^{-1}zw)^3=(x^{-1}z)^3w$.
But, since $(xz)^3=1$ and $z$ is of order $2$, we have $(x^{-1}z)^3=1$, which
implies $w=1$. This is absurd.

\medskip
{\sl Case} (b) {\sl with} $|G|=72$ ($|K|=3$):
Since $K$ is normal and $\Aut(\Z_3)\cong\Z_2$, we have $xwx^{-1}=w$, i.e, 
$xw=wx$.
We have $zwz=w$, or $w^{-1}$. If $zwz=w^{-1}$, then $w=(xz)^{-3}w(xz)^3=
zwz=w^{-1}$ (recall: $(xz)^3=1$), which is absurd.
Hence $zwz=w$.
Therefore, $K$ belongs to the center of $\langle x,z,w\rangle$.
Since $\rho(yzyz)=1$, we have $yzyz=1$, $w$, or $w^{-1}$. 
The case $yzyz=1$ is dismissed by the same reasoning as before.
If $yzyz=w$, then $w^2=z^2w^2=(zw)^2=(yzy)(yzy)=1$, contradiction.
The other case $yzyz=w^{-1}$ can be dismissed similarly. \qed

\vv

\vv

\paragraph\label{thm-groupthm2} {\bf Theorem.} \ {\sl If $G$ is a group of order 60 with more than one 5-Sylow group, then $G \cong A_5$.}

\vv

\pf The number of 5-Sylow groups of $G$ is then 6. The number of 3-Sylow groups of $G$ is either 1, 4 or 10. 

If there is a unique 3-Sylow, then the quotient by it is a group of order 20, hence has a unique 5-Sylow that pulls back to a normal subgroup of order 15 in $G$. This group has to be $\Z_{15}$, so its 5-Sylow is unique. Hence also $G$ has a unique 5-Sylow (here, we use that Sylows are conjugate, so if one of them belongs to a normal subgroup (here, $\Z_{15}$) of $G$, then all do).
If there are 4 3-Sylows, then we get a transitive action $\phi : G \rightarrow S_4$, of which the kernel is the intersection of the normalizers of these Sylows. The normalizers are of order 15, so this kernel is of order 5, 10 or 15. All groups of these orders have a unique 5-Sylow, which is then also normal in $G$ (by the same argument as before). 

Hence there are 10 3-Sylows. The number of 2-Sylows in $G$ is 1,3,5 or 15. If there is only one, the corresponding quotient is $\Z_{15}$, which has only one 5-Sylow, and the same argument shows this is impossible. If there are 3, then there is an action $G \rightarrow S_3$ with kernel of order 10 or 20, but both of these have a unique 5-Sylow, so the same applies. Now 15 2-Sylows don't fit
into $G$ together with the other Sylows, so $G$ has 5 2-Sylows. 

Hence there is an action $G \rightarrow S_5$; let $K$ be its kernel. The normalizer of a 2-Sylow has 12 elements, so $K$ has 1,2,3,4 or 6 elements, but if it has four, all normalizers have a common 2-Sylow (of order four). If it is 3 or 6, then $K$ contains a unique 3-Sylow, and the previous argument dismisses this possibility. If $K$ has order two, the image has order 30. This image has 1 or 6 5-Sylows, and 1 or 10 3-Sylows, so counting orders, there is at least a
normal 5-Sylow or a normal 3-Sylow. Pulling back to $G$ gives a normal 
subgroup of order 6 or 10, hence a unique 3- or 5-Sylow in $G$, a contradiction. 
Hence $K$ is trivial, and $G$ is a subgroup of $S_5$ of order 60, so $G \cong A_5$. \qed 

\vv

\noindent {\bf Acknowledgments.} The authors thank Prof.\ Hikoe Enomoto of Keio University for providing the proof of theorem \ref{thm-groupthm}. Part of this work was done at the MPIM (Bonn), and during a visit of the second author at Utrecht made possible in part by NWO.

\vv

\bibliographystyle{amsplain}
\providecommand{\bysame}{\leavevmode\hbox to3em{\hrulefill}\thinspace}

\vv

\vv

\noindent Dept.\ of Mathematics,
Utrecht University, P.O. Box 80010, 3508 TA  Utrecht, The Netherlands\footnote{All correspondence should be sent to this author}

\vv

\noindent Kyoto University, Faculty of Science, Dept.\ of Mathematics, Kyoto 606-8502, Japan

\vv

\noindent e-mail: {\tt cornelissen@math.uu.nl}, {\tt kato@kusm.kyoto-u.ac.jp}
}

\end{document}